# FIRST ORDER NORMALISATION IN THE GENERALISED PHOTOGRAVITATIONAL RESTRICTED THREE BODY PROBLEM WITH POYNTING-ROBERTSON DRAG

**B.S. Kushvah[1], J.P. Sharma[2] and B. Ishwar[3]**

[1] JRF DST Project, [2] Co-PI DST Project, [3] PI. DST Project

University Department of Mathematics,

B.R.A. Bihar University Muzaffarpur – 842001, India

Email: bskush@hotmail.com

In this paper we have studied non-linear stability of triangular equilibrium points. We have performed first order normalization in the generalized photogravitational restricted three body problem with Poynting-Robertson drag. In this problem we have taken bigger primary as source of radiation and smaller primary is an oblate spheroid. At first we have expanded the Lagrangian function in power series of x and y, where (x, y) are the co-ordinates of the triangular equilibrium points. Then the relation between the roots of the characteristic equation for the linearised system is obtained.

## 1. Introduction

The restricted three body problem describes the motion of an infinitesimal mass moving under the gravitational effect of the two finite masses, called primaries, which move in circular orbits around their centre of mass on account of their mutual attraction and the infinitesimal mass not influencing the motion of the primaries. The classical restricted three body problem is generalized to include the force of radiation pressure, the Poynting – Robertson effect and oblate ness effect.

J. H. Poynting (1903) considered the effect of the absorption and subsequent re-emission of sunlight by small isolated particles in the solar system. His work was later modified by H. P. Robertson (1937) who used a precise relativistic treatments of the first order in the ratio of he velocity of the particle to that light.

The effect of radiation pressure and P. R. drag in the restricted three body problem has been studied by Colombo *et at.* (1966), Chernikov Yu. A. (1970) and Schuerman (1980) who discussed the position as well as the stability of the Lagrangian equilibrium points when radiation pressure, P-R drag force are included. Murray C. D. (1994) systemematically discussed the dynamical effect of general drag In the planer circular restricted three body problem, Liou J. C. *et al.* (1995) examined the effect of radiation pressure, P-R drag and solar wind drag in the restricted three body problem.

Moser's conditions (1962), Arnold's theorem (1961) and Liapunov's theorem (1956) played a significant role in deciding the nonlinear stability of an equilibrium point. Applying Arnold's theorem (1961) Leontovic (1962) examined the nonlinear stability of triangular points. Moser gave some modifications in Arnold's theorem. Then Deprit and Deprit (1967) investigated the nonlinear stability of triangular points by applying Moser's modified version osf Arnold's theorem (1961).

Bhatnagar and Hallan (1983) studied the effect of perturbations on the nonlinear stability of triangular points. Maciejewski and Gozdziewski (1991) described the normalization algorithms of Hamiltonian near an equilibrium point. Further, Bhatnagar examined the nonlinear stability of $L_4$. Niedzielska (1994) investigated the nonlinear stability of the libration points in the photogravitational restricted three body problem. Mishra P and Bhola Ishwar (1995) studied the second order normalization in the generalized restricted problem of three bodies, smaller primary being an oblate spheroid.

In this paper we discuss on the first order normalization in the generalized photo gravitational restricted three body problem with Poynting-Robertson drag.

Further we will study non-linear stability of the triangular equilibrium point. For this we will apply Arnold's theorem (1961) and follow the procedure as adopted by Bhatnagar and Hallan (1983), Arnold proved that if :

(i) $k_1\omega_1 + k_2\omega_2 \neq 0$ for all pairs $(\kappa_1, \kappa_2)$ of rational integers and

(ii) Determinant $D \neq 0$, where $\omega_1, \omega_2$ are the basic frequencies for the linear dynamical system,

$$D = \det(bij), \quad (i,\ j = 1,\ 2,\ 3)$$

$$bij = \left(\frac{\partial^2 H}{\partial I_i \partial I_j}\right)_{Ii=Ij=0} \qquad (i,\ j = 1,\ 2)$$

$$b_{i3} = b_{3i} = \left(\frac{\partial H}{\partial I_i}\right)_{Ii=Ij=0} \qquad (i = 1,\ 2),$$

$$b_{33} = 0,$$

$$H = \omega_1 I_1 - \omega_2 I_2 + \frac{1}{2}\left[A I_1^2 + 2B I_1 I_2 + C {I_2}^2\right] + .......$$

is the normalized Hamiltonian with $I_1$ and $I_2$ as the action momenta coordinates, then an each energy manifold $H = h$ in the neighbourhood of equilibrium, there exist invariant tori of quasi-periodic motion, which divide the manifold, and consequently the equilibrium is stable. This is valid for a system with two degrees of freedom, which is the case under consideration. Moser has shown that Arnold's theorem is true if the condition (i) of the theorem is replaced by $k_1\omega_2 + k_2\omega_2 \neq 0$ for all pairs $(k_1, k_2)$ rational integers such that $|k_1| + |k_2| \leq 4$.

## 2. First Order Normalization

Equations of motion are

$$U_x = \ddot{x} - 2n\dot{y} = n^2 x - \frac{(1-\mu)q_1(x+\mu)}{r_1^3} - \frac{\mu(x+\mu-1)}{r_2^3} - \frac{3}{2}\mu A_2 \frac{(x+\mu-1)}{r_2^5}$$

$$-\frac{W_1}{r_1^2}\left\{\frac{(x+\mu)}{r_1^2}\left[(x+\mu)\dot{x} + y\dot{y}\right] + \dot{x} - ny\right\} \qquad \ldots (1)$$

$$U_y = \ddot{y} + 2n\dot{x} = n^2 y - \frac{(1-\mu)q_1 y}{r_1^3} - \frac{\mu y}{r_2^3} - \frac{3}{2}\frac{\mu A_2 y}{r_2^5}$$

$$-\frac{W_1}{r_1^2}\left\{\frac{y}{r_1^2}\left[(x+\mu)\dot{x} + y\dot{y}\right] + \dot{y} + n(x+\mu)\right\} \qquad \ldots (2)$$

where $W_1 = \dfrac{(1-\mu)(1-q_1)}{Cd}, \quad n^2 = 1 + \dfrac{3}{2}A_2, \; q = \left(1 - \dfrac{Fp}{Fg}\right),$

$$q = 1 - \frac{5.6\times10^{-5}}{as}\chi, \quad r_1^2 = (x+\mu)^2 + y^2, \quad r_2^2 = (x+\mu-1)^2 + y^2,$$

$$\mu = \frac{m_2}{m_1 + m_2} \leq \frac{1}{2}.$$

Where $m_1, m_2$ are masses of the primaries. The perturbed mean motion of the primaries is *n* and force function is

$$U = \frac{n^2}{2}\left(x^2 + y^2\right) + \frac{(1-\mu)q_1}{r_1} + \frac{\mu}{r_2} + \frac{1}{2}\frac{\mu A_2}{r_2^3}$$

$$+W_1\left\{\frac{(x+\mu)\dot{x} + y\dot{y}}{2r_1^2} - n\, arc\tan\left(\frac{y}{x+\mu}\right)\right\} \qquad \ldots (3)$$

The Lagrangian function of the problem can be written as

$$L=\frac{1}{2}\left(\dot{x}^2+\dot{y}^2\right)+n\left(x\dot{y}-\dot{x}y\right)+\frac{n^2}{2}\left(x^2+y^2\right)+\frac{(1-\mu)q_1}{r_1}+\frac{\mu}{r_2}$$

$$+\frac{\mu A_2}{2r_2^3}+W_1\left\{\frac{(x+\mu)\dot{x}+y\dot{y}}{2r_1^2}-n\,arc\tan\left(\frac{y}{x+\mu}\right)\right\}, \qquad \ldots (4)$$

and the Hamiltonian H is

$$H=-L+p_x\dot{x}+p_y\dot{y}=\frac{1}{2}\left(\dot{x}^2+\dot{y}^2\right)-\frac{n^2}{2}\left(x^2+y^2\right)-\frac{(1-\mu)q_1}{r_1}-\frac{\mu}{r_2}$$

$$-\frac{\mu A_2}{2r_2^3}+nW_1\left\{arc\tan\frac{y}{x+\mu}\right\}, \qquad \ldots (5)$$

where $P_x$, $P_y$ are the momenta coordinates given by

$$\left.\begin{aligned} p_x&=\frac{\partial L}{\partial \dot{x}}=\dot{x}-ny+\frac{W_1}{2r_1^2}(x+\mu)\\ p_y&=\frac{\partial L}{\partial \dot{y}}=\dot{y}+nx+\frac{W_1}{2r_1^2}y \end{aligned}\right\} \qquad \ldots (6)$$

The coordinates of the triangular equilibrium points are

$$x_*=x_0\left\{1-\frac{nW_1\left[(1-\mu)\left(1-\frac{5}{2}A_2\right)+\mu\left(1-\frac{A_2}{2}\right)\frac{\delta^2}{2}\right]}{3\mu(1-\mu)y_0x_0}-\frac{\delta^2A_2}{2x_0}\right\}$$

$$y_*=y_0\left\{1-\frac{nW_1\delta^2\left[2\mu-1-\mu\left(1-\frac{3}{2}A_2\right)\frac{\delta^2}{2}+7(1-\mu)\frac{A_2}{2}\right]}{3\mu(1-\mu)y_0^3}-\frac{\delta^2\left(1-\frac{\delta^2}{2}\right)A_2}{y_0^2}\right\}^{1/2}$$

where $x_0=\frac{\delta^2}{2}-\mu,\quad y_0=\pm\delta\left(1-\frac{\delta^2}{4}\right)^{1/2},\quad \delta=q^{1/3}$ … (7)

We shift the origin to $L_4$ for that, we change

$$x\to x_*+x,\quad y\to y_*+y,$$

The Lagrangian function L becomes

$$L=\frac{1}{2}\left(\dot{x}^2+\dot{y}^2\right)+n\left\{(x_*+x)\dot{y}-\dot{x}(y_*+y)\right\}+\frac{n^2}{2}\left\{(x_*+x)^2+(y_*+y)^2\right\}$$

$$+\frac{(1-\mu)q_1}{r_1}+\frac{\mu}{r_2}+\frac{\mu A_2}{2r_2^3}+W_1\left\{\frac{(x_*+x+\mu)\dot{x}+(y_*+y)\dot{y}}{2r_1^2}-n\;arc\tan\frac{(y_*+y)}{(x_*+x+\mu)}\right\}$$

… (8)

Where

$$r_1^{-1} = \left[(x+a)^2 + (y+b)^2\right]^{1/2} = f_1(x,\, y) \;\{\text{say}\}$$

$$r_2^{-1} = \left[(x+a-1)^2 + (y+b)^2\right]^{-1/2} = f_2(x,\, y)$$

$$r_1^{-2} = \left[(x+a)^2 + (y+b)^2\right]^{-1} = f_3(x,\, y)$$

$$r_2^{-3} = \left[(x+a-1)^2 + (y+b)^2\right]^{-3/2} = f_4(x,\, y)$$

$$arc\tan\left(\frac{y+b}{x+a}\right) = f_5(x,\, y)$$

$$a = x_* + \mu, \text{ and } \; b = y_*$$

Expanding L in power series of x and y, we get

$$\left.\begin{aligned} &L = L_0 + L_1 + L_2 + L_3 + L_4 + \ldots\ldots\ldots\ldots\ldots\ldots\ldots\ldots \\ &\text{or} \\ &H = H_0 + H_1 + H_2 + H_3 + H_4 + \ldots\ldots = -(L_0 + L_1 + L_2 \ldots.) + p_x\dot{x} + p_y\dot{y} \end{aligned}\right\}$$

… (9)

where

$$L_0 = \frac{n^2}{2}\left(a^2 + b^2\right) + (1-\mu)q_1\left(a^2 + b^2\right)^{-1/2} + \mu\left[(a-1)^2 + b^2\right]^{-1/2}$$

$$+\frac{\mu A_2}{2}\left[(a-1)^2 + b^2\right]^{-3/2} - nW_1 arc\tan\frac{b}{a}$$

$$L_1 = n(a\dot{y} - \dot{x}b) + n^2(ax+by) + x\{(1+\mu)q_1\left(-f_1^3 a\right) + \mu(a-1)f_2^3$$

$$+\frac{\mu A_2}{2}\left[-3f_2^5(a-1)\right] + nW_1 bf_3\} + y\{(1-\mu)q_1\left[-f_1^3.b\right] + \mu\left[-f_2^3 b\right]$$

$$+\frac{\mu A_2}{2}\left[-3f_2^5 b\right] + nW_1 f_3 a\} + W_1(a\dot{x} + b\dot{y})f_3$$

$$L_2 = \frac{1}{2}\left(\dot{x}^2 + \dot{y}^2\right) + n(x\dot{y} - \dot{x}y) + \frac{n^2}{2}\left(x^2 + y^2\right) + \frac{1}{\underline{|2}}\{x^2\left[(1+\mu)q_1\left(3f_1^2 a^2 - 1\right)f_1^3\right.$$

$$\left.+\mu\left[3f_2^3(a-1)^2 - 1\right]f_2^3 + \frac{\mu A_2}{2}\left[15f_2^7(a-1)^2 - 45f_2^5\right] - nW_1\left[2f_3^2 ab\right]\right]$$

$$+2xy\left[(1-\mu)q_1\left[6f_1^5ab\right]+\mu\left[6f_2^5(a-1)b\right]+\frac{\mu A_2}{2}\left[15f_2^7(a-1)b\right]\right.$$

$$-\frac{nW_1}{2}\left[2f_3^2b^2-2f_3^2a^2\right]\Big]+y^2\left[(1-\mu)q_1\left[3f_1^5b^2-f_1^3\right]+\mu\left[3f_2^5b^2-f_2^3\right]\right.$$

$$\left.+\frac{\mu A_2}{2}\left[15f_2^7b^2-45f_2^5\right]+nW_12f_3^2ab\right]\Big\}+W_1\left\{\frac{x\dot{x}+y\dot{y}}{2}f_3\right.$$

$$\left.-\left(a^2x\dot{x}+b^2y\dot{y}\right)f_3^3-\left[x\dot{y}ab+y\dot{x}ab\right]f_3^2\right\}$$

$$L_3=\frac{1}{L_3}\left\{x^3\left[(1-\mu)q_1\left(-15f_1^7a^3+9f_1^5a\right)+\mu\left[-15f_2^7(a-1)^3+9f_2^5(a-1)\right]\right.\right.$$

$$\left.+\frac{\mu A_2}{2}\left[-105f_2^9(a-1)^3+45(a-1)f_2^7\right]-nW_1\left[-8f_3^3a^2b^2+2f_3^2b\right]\right]$$

$$+3x^2y\left[(1-\mu)q_1\left(-15f_1^7a^2b+3f_1^5b\right)+\mu\left[-15f_2^7(a-1)^2+3f_2^5b\right]\right.$$

$$\left.+\frac{\mu A_2}{2}\left[-105f_2^9(a-1)^2b+15bf_2^7\right]-\frac{nW_1}{3}\left[-8f_3^3ab^2-2f_3^2a\right]\right]$$

$$+3xy^2\left[(1-\mu)q_1\left(-15f_1^7ab^2+3f_1^5a\right)+\mu\left(-15f_2^7b^2+f_2^5\right)(a-1)\right.$$

$$\left.+\frac{\mu A_2}{2}\left[-105f_2^9(a-1)b^2+15(a-1)f_2^7\right]-\frac{nW_1}{3}\left[-8f_3^3b^2+16f_3^3a^2b+2f_3^2b\right]\right]$$

$$+y^3\left[(1-\mu)q_1\left(-15f_1^7b^3+9f_1^5b\right)+\mu\left(-15f_2^7b^2+9f_2^5\right)b\right.$$

$$\left.\left.+\frac{\mu A_2}{2}\left[-105f_2^9b^3+45f_2^7b\right]-nW_1\left[8f_3^3ab^2-2f_3^2a\right]\right]\right\}$$

$$+W_1\frac{(a\dot{x}+b\dot{y})}{2}\left\{\frac{1}{\underline{|2}}\left[x^2\left(8f_3a^2-2\right)+16f_3abxy+\left[sf_3b^2-2\right]y^2\right]f_3^2\right\}$$

$$L_4=\frac{1}{L_4}\left\{x^4\left[\left(105f_1^9a^4-90f_1^7a^2+9f_1^5\right)(1-\mu)q_1+\mu\left[105f_2^9(a-1)^4-90f_2^7(a-1)^2+9f_2^5\right]\right.\right.$$

$$\left.+\frac{\mu A_2}{2}\left[945f_2^{11}(a-1)^2-630f_2^9(a-1)^2+45f_2^7\right]-nW_1\left[48f_3^4a^3b-24f_3^3ab\right]\right]$$

$$+4x^3y\left[(1-\mu)q_1\left(105f_1^9a^3b-45f_1^7ab\right)+\mu\left[105f_2^9(a-1)^3b-45f_2^7(a-1)b\right]\right.$$

$$\left.+\frac{\mu A_2}{2}\left[945f_2^{11}(a-1)^3b-315f_2^9(a-1)b\right]-\frac{nW_1}{4}\left[-48f_3^4a^4+144f_3^4a^2b^2+48f_3^3a^2-24f_3^2\right]\right]$$

$$+6x^2y^2\left[(1-\mu)q_1\left(105f_1^9a^2b^2-12f_1^5\right)+\mu\left[105f_2^9(a-1)^2b^2-12f_2^5\right]\right.$$

$$+\frac{\mu A_2}{2}\left[945 f_2^{11}(a-1)^2 b^2-105 f_2^9\left[(a-1)^2 b^2\right]+15 f_2^7\right]$$

$$-\frac{nW_1}{4}\left[144 f_3^4 ab^3-144 f_3^4 a^3 b-48 f_3^3 ab\right]\Big]$$

$$+4xy^3\left[(1-\mu)q_1\left(105 f_1^9 ab^3-45 f_1^7 ab\right)+\mu\left[105 f_2^9(a-1)b^3-45 f_2^7(a-1)b\right]\right.$$

$$+\frac{\mu A_2}{2}\left[945 f_2^{11}(a-1)b^3-315 f_2^9(a-1)b\right]-\frac{nW_1}{4}\left[48 f_3^4 b^4-48 f_3^3 b^2-144 f_3^4 a^2 b^2+24 f_3^2\right]\Big]$$

$$+y^4\left[(1-\mu)q_1\left(105 f_1^9 b^4-90 f_1^7 b^2+9 f_1^5\right)+\mu\left[105 f_2^9 b^4-90 f_2^7 b^2+9 f_2^5\right]\right.$$

$$+\frac{\mu A_2}{2}\left[945 f_2^{11} b^4-630 f_2^9 b^2+45 f_2^7\right]-\frac{nW_1}{4}\left[-48 f_3^4 ab^3+24 f_3^3 ab\right]\Big]\Big\}\Big\}$$

$$+W_1\frac{(a\dot{x}+b\dot{y})}{2.\underline{|3}}\left\{\left(-48 f_3 a^2+24\right)ax^3+3\left[-48 f_3 ab+8\right]bx^2 y\right.$$

$$\left.+3\left[-48 f_3 a+8\right]axy^2+\left[-48 f_3 b^2+24\right]b^3\right\}f_3^3$$

The second order part $H_2$ of the Hamiltonian H corresponding to Lagrangian L given in (8) takes the form

$$H_2=\frac{1}{2}\left(p_x^2+p_y^2\right)+n\left(yp_x-xp_y\right)+Ex^2+Fy^2+Gxy \qquad \ldots (10)$$

where $E=-\frac{1}{2}\left\{(1-\mu)q_1\left(2a^2-b^2\right)f_1^5+\mu\left[2(a-1)^2-b^2\right]f_2^5\right.$

$$\left.-\frac{15\,\mu A_2}{2}\left[\left(2(a-1)^2+3b^2\right)\right]f_2^7+2nW_1 f_3^2 ab\right\} \qquad \ldots (11)$$

$$F=-\frac{1}{2}\left\{(1-\mu)q_1\left(2b^2-a^2\right)f_1^5+\mu\left[2b^2-(a-1)^2\right]f_2^5\right.$$

$$\left.-\frac{15\mu A_2}{2}2\left[b^2+(a-1)^2\right]f_2^7+2nW_1 f_3^2 ab\right\} \qquad \ldots (12)$$

$$G=-\left\{6(1-\mu)q_1 f_1^5 ab+6\mu f_2^5(a-1)b\right.$$

$$\left.+15\mu A_2(a-1)bf_2^7-nW_1\left(b^2-a^2\right)f_3^2\right\} \qquad \ldots (13)$$

**3. Perturbed Basic Frequencies**

Now we follow the method of Whittaker (1965) to find the canonical transformation from the phase space $(x,\ y,\ p_x,\ p_y)$ in to the phase space product of the angle co-ordinates $(\phi_1, \phi_2)$ and the action momenta $(I_1, I_2)$ and of the first order in $I_1^{1/2}, I_2^{1/2}$, so that the second order part of the Hamiltonian is normalized. In order to examine the stability of the motion, we consider the following set of linear equations of the variables x, y :

$$-\lambda\ p_x = \frac{\partial H_2}{\partial x} = 2Ex + Gy - np_y$$

$$-\lambda\ p_y = \frac{\partial H_2}{\partial y} = 2Fy + Gx + np_x$$

$$\lambda_x = \frac{\partial H_2}{\partial p_x} = p_x + ny$$

$$\lambda_y = \frac{\partial H_2}{\partial p_y} = p_y - nx$$

i.e., $AX = 0$ … (14)

where, $A = \begin{bmatrix} 2E & G & \lambda & -n \\ G & 2F & n & \lambda \\ -\lambda & n & 1 & 0 \\ -n & -\lambda & 0 & 1 \end{bmatrix}$ and $X = \begin{bmatrix} x \\ y \\ p_x \\ p_y \end{bmatrix}$

since $X \neq 0 \Rightarrow |A| = 0$.

i.e., $\lambda^4 + 2\left(E + F + n^2\right)\lambda^2 + 4EF - G^2 + n^4 - 2n^2(E + F) = 0$ … (15)

This is the characteristic equation whose discriminant is

$$D = 4\left(E + F + n^2\right)^2 - 4\left\{4EF - G^2 + n^4 - 2n^2(E + F)\right\}$$

Stability is possible only when $D > 0$.We have seen, the Lagrangian function L and the Hamiltonian are the functions of *q,* $A_2$ and $W_1$. Hence they are affected by radiation pressure force, oblateness and P–R drag. So we conclude that these forces affect stability of triangular equilibrium points.

**Acknowledgement**

We are thankful to D.S.T, Govt. of India, New Delhi for sanctioning a project DST/MS/140/2K on this topic.